\documentclass[12pt]{article}
\usepackage{amsmath}
\usepackage{amssymb}
\usepackage{amsthm,amsxtra}
\usepackage{epsf}
\usepackage{eepic}
\usepackage{amsfonts, amscd, euscript,
graphics,psfrag}
\newlength{\mytopmargin}
\newlength{\myleftmargin}
\newtheorem{cor}{Corollary}

\newtheorem{prop}{Proposition}
\newcommand\psymmU{%
\begin{picture}(1,1)(0,0)%
\allinethickness{0.5pt}%
\path(0,0)(0,1)(1,1)(1,0)(0,0)%
\end{picture}}
\newcommand\psymmUU{%
\begin{picture}(1,1)(0,0)%
\allinethickness{0.5pt}%
\path(0,0)(0,1)(1,1)(1,0)(0,0)%
\put(0.5,0.5){\makebox(0,0){$\cdot$}}%
\end{picture}}
\newcommand\psymmO{%
\begin{picture}(1,1)(0,0)%
\allinethickness{0.5pt}%
\path(0,0)(0,1)(1,1)(1,0)(0,0)%
\path(0,0)(1,1)%
\end{picture}}
\newcommand\psymmS{%
\begin{picture}(1,1)(0,0)%
\allinethickness{0.5pt}%
\path(0,0)(0,1)(1,1)(1,0)(0,0)%
\path(1,0)(0,1)%
\end{picture}}
\newcommand\psymmu{%
\begin{picture}(1,1)(0,0)%
\allinethickness{0.5pt}%
\path(0,0)(0,1)(1,1)(1,0)(0,0)%
\path(0,0)(1,1)%
\path(0,1)(1,0)%
\end{picture}}

\newbox\tsymmUbox
\newbox\tsymmUUbox
\newbox\tsymmObox
\newbox\tsymmSbox
\newbox\tsymmubox
\setbox\tsymmUbox =\hbox{\kern0.75pt\setlength{\unitlength}{6pt}\psymmU \kern0.75pt}

\setbox\tsymmUUbox=\hbox{\kern0.75pt\setlength{\unitlength}{6pt}\psymmUU\kern0.75pt}
\setbox\tsymmObox =\hbox{\kern0.75pt\setlength{\unitlength}{6pt}\psymmO \kern0.75pt}
\setbox\tsymmSbox =\hbox{\kern0.75pt\setlength{\unitlength}{6pt}\psymmS \kern0.75pt}
\setbox\tsymmubox =\hbox{\kern0.75pt\setlength{\unitlength}{6pt}\psymmu \kern0.75pt}

\newbox\symmUbox
\newbox\symmUUbox
\newbox\symmObox
\newbox\symmSbox
\newbox\symmubox
\setbox\symmUbox =\hbox{\kern0.75pt\setlength{\unitlength}{4.5pt}\psymmU \kern0.75pt}
\setbox\symmUUbox=\hbox{\kern0.75pt\setlength{\unitlength}{4.5pt}\psymmUU\kern0.75pt}
\setbox\symmObox =\hbox{\kern0.75pt\setlength{\unitlength}{4.5pt}\psymmO \kern0.75pt}
\setbox\symmSbox =\hbox{\kern0.75pt\setlength{\unitlength}{4.5pt}\psymmS \kern0.75pt}
\setbox\symmubox =\hbox{\kern0.75pt\setlength{\unitlength}{4.5pt}\psymmu \kern0.75pt}

\def\symmUU{{\copy\symmUUbox}}
\def\symmO{{\copy\symmObox}}

\def\symmu{{\copy\symmubox}}

\newcommand{\Sc}{{\hbox{\rm Sc}}}
\newcommand{\Rc}{{\hbox{\rm Rc}}}
\newcommand{\HS}{{\hbox{\rm HS}}}
\newcommand{\E}{{\mathbb E}}

\newcommand{\Eu}{\E_{U_N}}
\newcommand{\natls}{{\mathbb N}}

\newcommand{\beq}{\begin{equation}}

\newcommand{\eeq}{\end{equation}}
\newcommand{\B}{\mathcal{B}}

\begin{document}
\vspace{4cm}
\noindent
{\bf Counting formulas associated with some random matrix averages}

\vspace{5mm} \noindent Peter J.~Forrester${}^{*}$ and Alex
Gamburd${}^\dagger$

\noindent
${}^*$Department of Mathematics and Statistics,
University of Melbourne, \\
Victoria 3010, Australia ; \\
${}^\dagger$ Department of Mathematics, University of California,
Santa Cruz \\and  Department of Mathematics, Stanford University,
\\Stanford, CA 94305 USA
%\email{agamburd@math.stanford.edu}\\

\small
\begin{quote}
Abstract: Moments of secular and inverse secular coefficients,
averaged over random matrices from classical groups, are related to
the enumeration of non-negative matrices with prescribed row and
column sums. Similar random matrix averages are related to certain
configurations of vicious random walkers and to the enumeration of
plane partitions. The combinatorial meaning of the average of the
characteristic polynomial of random Hermitian and Wishart matrices
is also investigated, and consequently several simple
universality results are derived.
\end{quote}

\section{Introduction}
The richness of random matrix theory was greatly enhanced during
the last couple of years of the millennium by the discovery of its
intimate connections with increasing subsequences and
non-intersecting lattice paths in enumerative and asymptotic
combinatorics (for reviews see \cite{AD99, Fo03}). These topics
relate to non-negative integer matrices via the celebrated
Robinson-Schensted-Knuth (RSK) correspondence \cite{Fu97, knuth}.
A special class of non-negative matrices are so called magic
squares
--- they have the property that the sum of the elements in any row
or column is equal to a prescribed positive integer $j$ say. Such
matrices are natural objects in RSK theory, so one might expect a
relationship between magic squares and random matrices. This is
indeed the case, and its development forms the main theme of the
recent work \cite{DiGa} (see \cite{BCAG} for related results). One
of the goals of the present paper is to continue the development
of this theme.

Magic squares are special cases of non-negative integer $m \times
n$ rectangular matrices in which the sum of the elements in each
row $j$ is equal to $\mu_j$, while the sum of the elements in each
column $k$ is equal to $\tilde{\mu}_k$. Without loss of
generality, for the purpose of enumeration we can insist that
$\mu_1 \ge \cdots \ge \mu_m$ and $\tilde{\mu}_1 \ge \cdots \ge
\tilde{\mu}_n$ so that $\mu := (\mu_1,\dots,\mu_m)$ and
$\tilde{\mu} := (\tilde{\mu}_1,\dots,\tilde{\mu}_n)$ form
partitions. We denote the total number of such matrices by $N_{\mu
\tilde{\mu}}$. Writing the partitions in terms of the frequencies
of their parts by
\begin{equation}\label{0.1}
\mu = \langle 1^{a_1} \cdots l^{a_l} \rangle, \qquad
\tilde{\mu} = \langle 1^{b_1} \cdots l^{b_l} \rangle
\end{equation}
where $l = {\rm max}(\mu_1,\tilde{\mu}_1)$, it was proved in
\cite[Thm.~2]{DiGa} that for $N \ge {\rm max}(\sum_{j=1}^l j a_j,
\sum_{j=1}^l j b_j)$
\begin{equation}\label{1}
N_{\mu \tilde{\mu}} =
{\mathbb E}_{M \in U(N)} \prod_{j=1}^l ( {\rm Sc}_j(M))^{a_j}
(\overline{{\rm Sc}_j(M))}^{b_j}.
\end{equation}
In (\ref{1}) the average is over matrices $M$ chosen from the group
$U(N)$ at random with respect to the Haar measure (uniform distribution), while
${\rm Sc}_j(M)$ is the $j$th secular coefficient of the characteristic
polynomial of $M$,
\begin{equation}\label{1.1}
P_M(z)=\det(M - zI_N) = \sum_{j=0}^N {\rm Sc}_j(M) (-z)^{N-j}.
\end{equation}
{}From (\ref{1}) it follows that the number $H_k(j)$ of $k \times k$
magic squares,
specified as $k \times k$
non-negative integer matrices with the sum of elements in
each row and column equal to $j$, is for $N \ge kj$, given by the formula
\begin{equation}
H_k(j) = {\mathbb E}_{M \in U(N)} | {\rm Sc}_j(M) |^{2k}.
\end{equation}

We will extend the formula (\ref{1}) in a number of directions.
The first relates to counting formulas analogous to (\ref{1}) in
the case that the matrices are constrained by a symmetry property.
Two such formulas, both relating to $n \times n$ symmetric
matrices, are known from \cite{DiGa}. Thus with the row sums
(which must be equal to the column sums as the matrices are
symmetric) labelled by the partition $\mu = (\mu_1,\dots,\mu_n) =
\langle 1^{a_1} \cdots l^{a_l} \rangle$ we have that for $|\mu| :=
\sum_{j=1}^l j a_j$ even, $N \ge |\mu|$, and all elements on the
diagonal zero, the total number $N_\mu^{\rm O}$ of such
non-negative integer matrices is given by
\begin{equation}\label{2}
N_{\mu}^{\rm O}
= {\mathbb E}_{M \in O(N)} \prod_{j=1}^l ( {\rm Sc}_j(M))^{a_j}.
\end{equation}
If instead all elements on the diagonal are permitted to be even, the total
number $N_\mu^{\rm Sp}$ of such matrices is given by
\begin{equation}\label{3}
N_{\mu}^{\rm Sp} = {\mathbb E}_{M \in USp(2N)} \prod_{j=1}^l (
{\rm Sc}_j(M))^{a_j}.
\end{equation}
In (\ref{2}) the average is over the group of orthogonal matrices
$O(N)$, while in (\ref{3}) it is over the group of $2N \times 2N$
unitary symplectic matrices $USp(2N)$ (both with the Haar
measure). In Section 2 we will extend (\ref{2}) and (\ref{3}) to
the case that the diagonal elements have no special restriction.
We also give analogous formulas in the case of square matrices
symmetric about both the diagonal and anti-diagonal, and for
matrices symmetric about the centre point of the matrix. That such
symmetrizations relate to averages over classical groups is a
consequence of results of Rains \cite{Ra98}, and  Baik and Rains
\cite{BR02}. Another tractable case in terms of a random matrix
average to be considered is the setting relating to (\ref{1}) but
with the entries of the matrix restricted to 0's and 1's.
If the matrix has a $2 \times 2$ block structure, with the two
diagonal blocks having non-negative integer entries, and the off
diagonal blocks having entries 0 or 1, counting formulas of the
type (\ref{1}) can be obtained by studying ratios of characteristic
polynomials averaged over the classical groups. This is done in
Section \ref{ratio}.

The topic of Section \ref{path} is the relationship between
enumeration formulas (\ref{1}) and (\ref{3}) and certain classes
of non-intersecting lattice paths, or equivalently certain
configurations of vicious random walkers. This is motivated by the
graphical representation of the RSK correspondence in terms of
non-intersecting lattice paths \cite{Jo02,FR02}.  In Section
\ref{part} we describe connections between moments of
characteristic polynomials and enumeration of certain classes of
plane partitions.  Finally, in Section \ref{herm}, we revisit the
question of the combinatorial meaning  of the expected value of
characteristic polynomials of random complex Hermitian Wigner
matrices, and take up a similar study in relation to Wishart
matrices.

\section{Further symmetrizations of the square}\label{s2}
\setcounter{equation}{0}
The RSK correspondence gives a bijection between weighted $n \times n$
non-negative integer matrices, entries $x_{i j}$ weighted by
$(\alpha_i \beta_j)^{x_{ij}}$, and pairs of weighted semi-standard
tableaux of content $n$. In the latter one member of the pair is
weighted by
\begin{equation}\label{4.1}
\alpha_1^{\# 1's} \alpha_2^{\# 2's} \cdots \alpha_n^{\# n's},
\end{equation}
while the other member is weighted by (\ref{4.1}) but with the
$\alpha$'s replaced by the $\beta$'s (the exponents in (\ref{4.1}) are
determined by the entries of the matrix). Summing (\ref{4.1}) over all
allowed semi-standard tableaux of a given shape $\kappa$ gives the
combinatorial definition of the Schur polynomial, and this way one
obtains for the generating function of the weighted matrices the
identity
\begin{equation}\label{4.2}
{1 \over \prod_{i,j = 1}^n (1 - \alpha_i \beta_j) } =
\sum_\kappa s_\kappa(\alpha_1,\dots,\alpha_n)
s_\kappa(\beta_1,\dots,\beta_n).
\end{equation}
This is well known in the theory of the Schur polynomial, and is called the
Cauchy formula (see e.g.~\cite{Mac}). We remark that
the case of $m \times n$ matrices ($m < n$ say for definiteness) follows
by simply setting $\alpha_{m+1}=\cdots =\alpha_n=0$.

The fact that the entries $x_{ij}$ are weighted by
$(\alpha_i \beta_j)^{x_{ij}}$ tells us
the coefficient of $$\alpha_1^{\eta_1} \cdots \alpha_n^{\eta_n}
\beta_1^{\rho_1} \cdots \beta_n^{\rho_n}$$ in (\ref{4.2}) counts
the number of matrices for which the sum of the elements in row
$i$ equals $\eta_i$, while the sum of the elements in row $j$
equals $\rho_j$. Because (\ref{4.2}) is symmetric in the
$\alpha_i$'s and the $\beta_j$'s, without loss of generality we
can restrict attention to the case that $\eta_1 \ge \cdots \ge
\eta_n$, $\rho_1 \ge \cdots \ge \rho_n$ and thus $\eta$ and $\rho$
form partitions, to be denoted $\mu$ and $\tilde{\mu}$
respectively say. The task is then to extract the coefficient of
$\alpha^\mu \beta^{\tilde{\mu}}$ (in an obvious multivariable
shorthand notation).

Only terms on the RHS of (\ref{4.2}) with $|\kappa| = |\mu| =
|\tilde{\mu}|$ can contribute to the coefficient of $\alpha^\mu
\beta^{\tilde{\mu}}$. Because of this, we can use a result of Baik and
Rains \cite{BR02} expressing the RHS of (\ref{4.2}) with the largest part
of $\kappa$ restricted to be no greater than $N$ as an average over
$U(N)$,
\begin{equation}\label{5.1}
\sum_{\kappa: \kappa_1 \le N}
s_\kappa(\alpha_1,\dots,\alpha_n) s_\kappa(\beta_1,\dots,\beta_n)
= {\mathbb E}_{U \in U(N)}
\prod_{j=1}^n \det(I_N + \alpha_j \bar{U})
\det(I_N + \beta_j U) .
\end{equation}
In particular with $\mu$, $\tilde{\mu}$ specified as in
(\ref{0.1}), (\ref{5.1}) agrees with (\ref{4.2}) at the required
order provided
\begin{equation}\label{2.4}
N \ge |\mu| = \sum_{j=1}^l j a_j = \sum_{j=1}^l j b_j.
\end{equation}
The extraction of the sought coefficient from (\ref{5.1}) is immediate,
and we reclaim (\ref{1}), provided $N$ obeys the inequality (\ref{2.4}).

We remark in passing, that  since we have \cite{CFKRS2, DiGa}
\beq \label{e:charmom1} {\mathbb E}_{U \in U(N)} \prod_{j=1}^m
\det(I_N + \alpha_j U) \prod_{j=1}^n \det(I_N + \overline{U}/\beta_j)
=\frac{1}{(\beta_1 \dots \beta_n)^N} s_{N^n}(\alpha_1, \dots,
\alpha_m; \beta_1, \dots, \beta_n), \eeq
(in fact this equation will be of independent interest below)
where ${N^n}$ denotes the partition with
$n$ parts all equal to $N$, it is also true that $N_{\mu \tilde{\mu}}$
is equal to the coefficient of $\alpha^\mu \beta^{N - \tilde{\mu}}$ in
the Schur function $S_{N^n}(\{\alpha_j,\beta_j\}_{j=1,\dots,n})$,
provided $N$ is large enough as required by (\ref{2.4}).

{}From the present perspective, to extend the enumeration formula
(\ref{1}) to classes of symmetrized non-negative matrices we
require a formula analogous to (\ref{5.1}) for the corresponding
partial generating function. Indeed such a formula (due to
Littlewood \cite{Li}) is known for symmetric matrices with all
elements zero on the diagonal,
\beq \label{lit1} \sum_{\lambda' \,
\text{even}} s_{\lambda}(\alpha_1, \dots ,\alpha_n)= \prod_{1 \le i
< j \le n}\frac{1}{1-\alpha_i \alpha_j}, \eeq
and  for symmetric
matrices with all elements even on the diagonal,
\beq \label{lit2}
\sum_{\lambda \, \text{even}} s_{\lambda}(\alpha_1, \dots ,
\alpha_n)= \prod_{1 \le i \le j \le n}\frac{1}{1-\alpha_i
\alpha_j}.\eeq
Combining \eqref{lit1} with
\beq \label{lo} {\mathbb E}_{M \in O(N)}
\prod_{j=1}^n \det (I_N + \alpha_j M) = \sum_{
\substack{l(\lambda) \le N \\ \lambda' \,
\text{even}}}s_{\lambda}(\alpha_1, \dots ,\alpha_n)
\eeq
and \eqref{lit2} with
\beq \label{lsp}
{\mathbb E}_{M \in USp(2N)} \prod_{j=1}^n \det (I_{2N} + \alpha_j M) =
\sum_{ \substack{\lambda_1 \le 2N \\ \lambda \,
\text{even}}}s_{\lambda}(\alpha_1, \dots , \alpha_n) \eeq
allows us
to reclaim (\ref{2}) and (\ref{3}) respectively. The finitizations
(\ref{lo}) and (\ref{lsp}) of Littlewood's formula are due to Rains
\cite{Ra98}.

 Formulas
analogous to (\ref{5.1}) are also known in the case of symmetric
matrices with no special constraint on the diagonal, matrices
symmetric about both the diagonal and anti-diagonal, and matrices
with a point reflection symmetry about the centre. Let us first
consider symmetric matrices with no special constraint on the
diagonal.

For a symmetric matrix, the row sum equals the column sum so we can take
for the weights
$\alpha_i \beta_j = \sqrt{q_i q_j}$ and extract the
coefficient of $q^\mu$ to read off the number of non-negative
integer matrices with
row sums equal to $\mu$. Since each off diagonal element in position
$ij$ of a
symmetric matrix is identical to the element in position
$ji$, we can equivalently
restrict attention to the lower triangular portion $i \le j$ of the matrix,
and weight the elements in
strictly lower triangular positions $i<j$ by $q_i q_j$, and
the diagonal elements by $q_i$. The generating function for weighted
symmetric matrices is thus seen to be
\begin{equation}\label{7.1}
{1 \over \prod_{i=1}^n (1 - q_i) \prod_{i < j} (1 - q_i q_j) }
\end{equation}
and we seek the coefficient of $q^\mu$ in this expression.
According to an identity of Littlewood \cite{Li} we have that
(\ref{7.1}) is equal to
\begin{equation}\label{7.1a}
\sum_{\kappa} s_\kappa(q_1,\dots,q_n).
\end{equation}
It has been shown by Baik and Rains \cite{BR02} that (\ref{7.1a})
constrained by the size of the largest part of $\kappa_1$ can be expressed
as a random matrix average,
\begin{equation}\label{7.1b}
\sum_{\kappa:\kappa_1 \le N } s_\kappa(q_1,\dots,q_n) =
{\mathbb E}_{U \in O(N)} \det(I_N + U)
\prod_{j=1}^n \det ( I_N + q_j U).
\end{equation}
With knowledge of the above results, we can now easily obtain an
enumeration formula for symmetric non-negative integer matrices with
prescribed row sums.

\begin{prop}\label{p1}
Consider the set of all symmetric $n \times n$ non-negative integer matrices.
Label the row sums (or equivalently column sums) by the partition
$$
\mu = (\mu_1,\dots,\mu_n) = \langle 1^{a_1} \cdots l^{a_l} \rangle.
$$
For $N \ge |\mu|$ the total number $N_\mu^{\symmO}$ of such matrices is
given as a random matrix average by
\begin{equation}\label{8.1}
N_\mu^{\symmO} = {\mathbb E}_{M \in O(N)} \det (I_N + M)
\prod_{j=1}^l ({\rm Sc}_j(M))^{a_j}
\end{equation}
(cf.~(\ref{2})).
\end{prop}

\noindent
Proof. \quad Only terms in (\ref{7.1a}) with $|\kappa| = |\mu|$ can contribute
to the coefficient of $q^\mu$ in (\ref{7.1}), so for $N \ge |\mu|$ the
latter is the same as the coefficient of $q^\mu$ in (\ref{7.1b}).
Extracting the coefficient gives (\ref{8.1}).
\hfill $\square$

We remark that since $\det(I_N+M)=0$ for $M \in O^-(N)$ (because
$\lambda = -1$ is an eigenvalue) we can restrict the matrices $M$ in
(\ref{8.1}) to $M \in O^+(N)$ (the matrices $U$ in (\ref{7.1b}) can similarly
be restricted).

\medskip
Next, we note that the results of Baik and Rains in fact allow the
diagonal sum to be prescribed. The relevant generating function is
\cite{BR02}
\begin{equation}
{1 \over \prod_{i=1}^n(1 - \alpha q_i) \prod_{i<j}(1 - q_i q_j) }
= \sum_{\kappa} \alpha^{\sum_{j=1}^n(-1)^{j-1} \kappa_j}
s_\kappa(q_1,\dots,q_n)
\end{equation}
and for $\kappa_1 \le N$ this sum can be written as the random matrix
average
\begin{equation}\label{2.16}
{\mathbb E}_{U \in O(N)} \det(I_N + \alpha U)
\prod_{j=1}^n \det ( I_N + q_j U).
\end{equation}
Hence, with the diagonal sum prescribed to be equal to $p$ say and
the total number of symmetric non-negative matrices now denoted
$N_{\mu,p}^\symmO$, the enumeration formula (\ref{8.1}) should be
modified to read
\begin{equation}
N_{\mu,p}^\symmO =
{\mathbb E}_{M \in O(N)} {\rm Sc}_p(M))
\prod_{j=1}^l ({\rm Sc}_j(M))^{a_j}.
\end{equation}

We now turn our attention to the task of enumerating according to row
sums $2n \times 2n$ non-negative integer matrices $[x_{i,j}]_{i,j=1,\dots,
2n}$ with a reflection symmetry about the diagonal, $x_{i,j} =
x_{j,i}$ $(i>j)$, and about the anti-diagonal,
$x_{i,j} = x_{i,2n+1-j}$ $(i>2n+1-j)$. The generating function for
weighted symmetric matrices of this type is
\begin{equation}\label{8.2}
\prod_{i=1}^n {(1+\chi_1 q_i) \over (1-\chi_0 q_i) }
{1 \over \prod_{i,j=1}^n(1 - q_i q_j) }
\end{equation}
where $\chi_0=1$ ($\chi_1=1$) in the case that the elements on the
diagonal (anti-diagonal) are unrestricted, while $\chi_0=0$
($\chi_1=0$) in the case the elements on the diagonal are restricted to be
zero (anti-diagonal are restricted to be even). We know from
\cite{BR02} that the random matrix average
\begin{equation}\label{8.3}
{\mathbb E}_{U \in U(N)}
{\det( 1 + \chi_0 U) \over \det( I_N - \chi_1 U) }
\prod_{j=1}^n \det (I_N + q_j U) \det(I_N + q_j \bar{U})
\end{equation}
correctly reproduces (\ref{8.2}) up to and including terms
$O(q_j^{2N+1})$ in the $q_j$. (In the case $\chi_1=1$, (\ref{8.3}) is
to be interpreted as the limiting value for $\chi_1 \to 1^-$.)
The following enumeration result is now evident.

\begin{prop}
Consider the set of all symmetric $2n \times 2n$ non-negative integer matrices
which have the further constraint of being symmetric about the anti-diagonal.
Introduce possible constraints on the diagonal and
anti-diagonal elements according to the values of
$\chi_0$ and $\chi_1$ noted below (\ref{8.2}). With the row sums labelled
by $\mu$, we have that for $N \ge {1 \over 2} |\mu|$ the number of
matrices in the set, $N_\mu^\symmu$ say, is given by the coefficient
of $q^\mu$ in the random matrix average (\ref{8.3}).
\end{prop}

As our next example of a symmetry constraint, we turn our
attention to the case of $2n \times 2n$ matrices invariant with
respect to reflections in the point $(n+1/2,n+1/2)$ (here we are
thinking of the matrix as labelled by a grid of $2n \times 2n$
lattice points $\{(i,j): 1 \le i,j \le 2n \}$). For such matrices
$x_{i,j} = x_{2n+1-i, 2n+1-j}$. The generating function for this
class of matrices is
\begin{equation}\label{9.1}
{1 \over \prod_{i,j=1}^n (1 - q_i q_j)^2 }.
\end{equation}
It follows from (\ref{4.2}) and (\ref{5.1}) that up and including terms
$O(q_i^{N})$ (\ref{9.1}) is equal to
\begin{equation}\label{10.1}
\Big ( {\mathbb E}_{U \in U(N)}
\prod_{j=1}^n \det (I_N + q_j \bar{U}) \det (I_N + q_j U) \Big )^2.
\end{equation}
We thus have the following enumeration result.

\begin{prop}
Consider the set of all $2n \times 2n$ non-negative integer matrices with
the point reflection symmetry $x_{i,j} = x_{2n+1-i, 2n+1-j}$. With the
row sums labelled by $\mu$, we have that for $N \ge
|\mu|$, the number of matrices in the set $N_\mu^\symmUU$ say is given
by the coefficient of $q^\mu$ in the random matrix average
(\ref{10.1}).
\end{prop}

As another extension of (\ref{1}) we consider not a symmetry
constraint on the matrix, but rather a restriction on the entries
of the matrix. These we take to be either $0$ or 1. For $m \times
n$ matrices of this type, weighted by $(\alpha_i
\beta_j)^{x_{ij}}$, the generating function is
\begin{equation}\label{2.21a}
\prod_{i=1}^m \prod_{j=1}^n (1 + \alpha_i \beta_j).
\end{equation}
According to the dual Cauchy identity,
\begin{equation}\label{dC}
\prod_{i,j=1}^n(1+\alpha_i \beta_j) =
\sum_\mu s_{\mu'}(\alpha_1,\dots,\alpha_n) s_\mu(\beta_1,\dots,\beta_n),
\end{equation}
where $\mu'$ denotes the partition conjugate to $\mu$ (see
e.g.~\cite{Mac}). Analogous to (\ref{5.1}), we have that this sum
restricted by the size of the largest part of $\mu$ can be written
as a random matrix average
\begin{equation}\label{dC1}
\sum_{\mu: \mu_1 \le N}
s_{\mu'}(\alpha_1,\dots,\alpha_n) s_\mu(\beta_1,\dots,\beta_n)
=
 {\mathbb E}_{U \in U(N)}
\prod_{j=1}^n {\det (I_N + \alpha_j {U}) \over
\det (I_N - \beta_j \bar{U}) }.
\end{equation}
Because of the occurrence  of the dual partition in (\ref{dC}),
$\mu_1$ cannot exceed $n$ for a non-zero contribution and so
(\ref{dC1}) is independent of $N$ for $N \ge n$.

We want to extract from this the coefficient of $\alpha^\mu
\beta^{\tilde{\mu}}$. For this purpose we introduce
inve{r}se secular coefficients $\Rc_p(U)$ according to
$$
{1 \over \det (I_N - xU) } = \sum_{p=0}^\infty x^p \Rc_p(U).
$$
In terms of the eigenvalues $e^{i\theta_j}$ $(j=1,\dots,N)$ of $U$ we have
that
$$
\Rc_p(U) = h_p(e^{i\theta_1},\dots,e^{i\theta_N})
$$
where $h_p$ denotes the $p$th complete symmetric function (see
\cite{Mac}). Making use too of (\ref{1.1}) we obtain for the
sought enumeration the following result.

\begin{prop}
Consider the set of all $m \times n$ matrices in which the entries
take on the values  0  or 1. Prescribe the row sums and column
sums by the partitions $\mu$ and $\tilde{\mu}$ respectively, which
are to be written in terms of the frequency of their parts by
(\ref{0.1}). For $N \ge {\rm min}(m,n,\sum_{j=1}^l j a_j,
\sum_{j=1}^l j b_j)$ the total number $N_{\mu\tilde{\mu}}^{0,1}$
of such matrices is given in terms of a random matrix average by
\begin{equation}\label{sr}
N_{\mu\tilde{\mu}}^{0,1} = {\mathbb E}_{U \in U(N)} \prod_{j=1}^l
({\rm Sc}_j(U))^{a_j}( \Rc_j(\bar{U}))^{b_j}.
\end{equation}
\end{prop}

\noindent
We remark that for (\ref{sr}) to be non-zero we require $\mu_1 \le m$,
$\tilde{\mu}_1 \le n$ and $|\mu| = |\tilde{\mu}|$.

\section{Block matrix structures and
ratios of characteristic polynomials} \label{ratio}
\setcounter{equation}{0}

We have seen in the previous section how generating functions
for non-negative integer matrices, and 0-1 matrices, are related
to Schur function identities which in turn are related to averages
over the classical groups. Results obtained in \cite{BR85,BR02}
tell us this strategy can also be carried through for classes of
block matrices which relate to averages over the unitary, symplectic
and orthogonal groups.

In relation to the unitary group, let the block structured
$(k_1+l_1) \times (k_2 + l_2)$ non-negative integer matrix
\begin{equation}\label{abcd}
\begin{bmatrix} A & B \\ C & D \end{bmatrix}
\end{equation}
be such that $A$ is a $k_1 \times k_2$ non-negative integer matrix;
$D$ is a $l_1 \times l_2$ non-negative integer matrix; $B$ is a
$k_1 \times l_2$ 0-1 matrix; and $C$ is a $l_1 \times k_2$ 0-1
matrix. Suppose the entries of the block $X = [x_{ij}]$ are weighted by
$g_{ij}^{x_{ij}}$ with
\begin{equation}\label{g}
g_{ij} = \left \{ \begin{array}{ll} \alpha_i \gamma_j & {\rm for}
\quad X = A \\
\beta_i \delta_j & {\rm for} \quad X = D \\
\alpha_i \delta_j & {\rm for} \quad X = B \\
\beta_i \gamma_j & {\rm for} \quad X = C.
\end{array} \right.
\end{equation}
Such matrices have for their generating function
$$
\prod_{i=1}^{k_1}\prod_{j=1}^{k_2}\frac{1}{1-\alpha_i \gamma_j}
\prod_{i=1}^{l_1}\prod_{j=1}^{l_2}\frac{1}{1-\beta_i \delta_j}
\prod_{i=1}^{k_1}\prod_{j=1}^{l_2}(1+\alpha_i \delta_j)
\prod_{i=1}^{l_1}\prod_{j=1}^{k_2}(1+\beta_i \gamma_j),
$$
thus containing as special cases both the LHS of (\ref{4.2}) and
(\ref{2.21a}). The coefficient of $\alpha^{\mu} \beta^{\tilde{\mu}}
\gamma^\nu \delta^{\tilde{\nu}}$ tells us the number of matrices
(\ref{abcd}) with prescribed row and column sums.
Moreover generalizations of both the Cauchy identities
(\ref{4.2}), (\ref{dC}), and their finitizations (\ref{5.1}),
(\ref{dC1}) are known \cite{BR85,BR02}, allowing for the
counting function to be expressed as a random matrix average.

To state these generalizations requires introducing the functions
\cite{Mac}
\beq
\HS_{\lambda}(\alpha_1, \dots, \alpha_{k};\beta_1,
\dots \beta_{l}) =\det(a_{\lambda_{i}+j-1})_{1\le i, j \le
l(\lambda)},
\eeq
where $a_k$ denotes the coefficient of $x^k$ in
$$
{ \prod_{j=1}^{l}(1+\beta_j
x) \over \prod_{i=1}^{k}(1-\alpha_i x) }.
$$
The HS${}_\lambda$ are referred to as the hook Schur functions. In terms of
this the generalization of the Cauchy identities is
\begin{eqnarray}\label{e:hs2}
&&\sum_{\lambda}\HS_{\lambda}(\alpha_1, \dots,
\alpha_{k_1};\beta_1, \dots \beta_{l_1}) \HS_{\lambda}(\gamma_1,
\dots, \gamma_{k_2};\delta_1, \dots \delta_{l_2})\nonumber \\&&
 \qquad = \prod_{i=1}^{k_1}\prod_{j=1}^{k_2}\frac{1}{1-\alpha_i \gamma_j}
\prod_{i=1}^{l_1}\prod_{j=1}^{l_2}\frac{1}{1-\beta_i \delta_j}
\prod_{i=1}^{k_1}\prod_{j=1}^{l_2}(1+\alpha_i \delta_j)
\prod_{i=1}^{l_1}\prod_{j=1}^{k_2}(1+\beta_i \gamma_j)
\end{eqnarray}
while the generalizations of their finitizations is
\begin{eqnarray} \label{e:hs1}
&&{\mathbb E}_{U \in U(N)} \frac {\prod_{i=1}^{k_1}
\det (I_N + \alpha_i U) \prod_{j=1}^{k_2} \det (I_N +\gamma_j
\bar{U})} {\prod_{m=1}^{l_1} \det (I_N -\beta_m U)
\prod_{n=1}^{l_2} \det (I_N -\delta_n \bar{U})} \nonumber \\&&
\qquad = \sum_{\lambda_1 \le N}\HS_{\lambda}(\alpha_1, \dots,
\alpha_{k_1};\beta_1, \dots \beta_{l_1}) \HS_{\lambda}(\gamma_1,
\dots, \gamma_{k_2};\delta_1, \dots \delta_{l_2}).
\end{eqnarray}

Arguing as in the derivation of the counting formula (\ref{1})
given in the first three paragraphs of Section \ref{s2}, the
following generalization of (\ref{1}) and (\ref{sr}) is immediate.

\begin{prop}
Let $\mathbf{a}=(a_1, \dots, a_{k_1})$ , $\mathbf{b}=(b_1,
\dots, b_{l_1})$,  $\mathbf{c}=(c_1, \dots, c_{k_2})$ ,
$\mathbf{d}=(d_1, \dots, d_{l_2})$,  where  the $a_j$, $b_j$, $c_j$,
$d_j$ are nonnegative integers.
{}From these arrays form partitions
$$
\mu=\langle 1^{a_1}\cdots
k_1^{a_{k_1}}\rangle , \: \: \tilde{\mu}=\langle 1^{b_1}\cdots
l_1^{b_{l_1}}\rangle, \: \:
\nu=\langle 1^{c_1}\cdots k_2^{c_{k_2}}\rangle, \: \:
\tilde{\nu}=\langle 1^{d_1}\cdots l_2^{d_{l_2}}\rangle.
$$
Let $N_{\mu \tilde{\mu} \nu \tilde{\nu}}$ denote the number of
matrices (\ref{abcd}) with
$$
{\rm row}(A,B) = \mu, \: \: {\rm row}(C,D) = \nu, \: \:
{\rm col}(A,C) = \tilde{\mu}, \: \: {\rm col}(B,D) = \tilde{\nu}
$$
where the notation row$(X,Y)$ refers to the row sums across $X$
and $Y$ given they are horizontal neighbors in a block matrix, and
col$(X,Y)$ refers to the column sums down $X$ and $Y$ given they
vertical neighbors in a block matrix.

For
$$
N \geq
\max\left(\sum_1^{k_1} ja_j, \sum_1^{l_1} j b_j, \sum_1^{k_2}
jac_j, \sum_1^{l_2} j d_j, \right)
$$
we have
\beq
\label{e:mixedmom} \Eu \prod_{i=1}^{k_1} (\Sc_i(M))^{a_i}
\prod_{j=1}^{l_1}\overline{(\Sc_j(M))} ^{b_j} \prod_{m=1}^{k_2}
(\Rc_m(M))^{c_m} \prod_{n=1}^{l_2}\overline{(\Rc_n(M))} ^{d_n} =
N_{\mu \tilde{\mu} \nu \tilde{\nu}}.
\eeq
\end{prop}

As reviewed in the Introduction, a relationship with averages over
the orthogonal and symplectic groups comes about when the matrix
is constrained to be symmetric. Thus in (\ref{abcd}) we must take
\begin{equation}\label{CB}
k_1 = k_2 = k, \quad l_1 = l_2 = l, \quad
A = A^T, \quad D = D^T, \quad C = B^T.
\end{equation}
Correspondingly, for $i \ne j$, the weights (\ref{g}) are to be
replaced by
\begin{equation}\label{gg}
g_{ij} = \left \{ \begin{array}{ll} \sqrt{ \alpha_i \alpha_j} & {\rm for}
\quad X = A \\
\sqrt{\beta_i \beta_j} & {\rm for} \quad X = D \\
\sqrt{\alpha_i \beta_j} & {\rm for} \quad X = B \\
\sqrt{\beta_i \alpha_j} & {\rm for} \quad X = C.
\end{array} \right.
\end{equation}
As in the meaning of (\ref{2}) and (\ref{3}), the constraint on
the diagonal elements will determine whether the relationship is
with an average over the symplectic group, or an average over the
orthogonal group. It turns out that the constraint relevant to the
symplectic (orthogonal) group is that all elements on the diagonal
of $A$ be even (zero), while those on $D$ be zero (even). The
corresponding generating functions are therefore given by
\begin{eqnarray*}
G^{\rm Sp}(\{\alpha_i\},\{\beta_j\}) & = &
\prod_{1 \le i \le j \le k}\frac{1}{1-\alpha_i \alpha_j}
 \prod_{1 \le i < j \le l}\frac{1}{1-\beta_i \beta_j}
\prod_{i=1}^{k}\prod_{j=i}^{l} (1+\alpha_i \beta_j) \\
G^{\rm O}(\{\alpha_i\},\{\beta_j\}) & = &
\prod_{1 \le i < j \le k}\frac{1}{1-\alpha_i \alpha_j}
 \prod_{1 \le i \le j \le l}\frac{1}{1-\beta_i \beta_j}
\prod_{i=1}^{k}\prod_{j=i}^{l} (1+\alpha_i \beta_j)
\end{eqnarray*}

In relation to these generating functions, one has as generalizations
of the Littlewood identities (\ref{lit1}) and (\ref{lit2})
\cite{BR85}
\begin{eqnarray*}
\sum_{\lambda \, \text{even}}
\HS_{\lambda}(\alpha_1, \dots \alpha_k; \beta_1, \dots, \beta_l)
& = & G^{\rm Sp}(\{\alpha_i\},\{\beta_j\}) \nonumber \\
\sum_{\lambda' \, \text{even}}
\HS_{\lambda}(\alpha_1, \dots \alpha_k; \beta_1, \dots, \beta_l)
& = & G^{\rm O}(\{\alpha_i\},\{\beta_j\}).
\end{eqnarray*}
And with the largest part of the partitions in the sum restricted, the
LHS's have the finitizations \cite{BR02}
\begin{eqnarray}  \label{e:hsp1} {\mathbb E}_{U \in USp(2N)} \frac
{\prod_{i=1}^{k} \det (I_{2N} + \alpha_i U)} {\prod_{j=1}^{l}
\det (I_{2N} -\beta_j U)} & = & \sum_{ \substack{\lambda_1 \le 2N \\
\lambda \, \text{even}}} \HS_{\lambda}(\alpha_1, \dots \alpha_k;
\beta_1, \dots, \beta_l) \\
\label{e:hso1} {\mathbb E}_{U \in O(N)} \frac
{\prod_{i=1}^{k} \det (I_N + \alpha_i U)} {\prod_{j=1}^{l}
\det (I_N -\beta_j U)} & = & \sum_{ \substack{l(\lambda) \le N \\
\lambda' \, \text{even}}} \HS_{\lambda}(\alpha_1, \dots \alpha_k;
\beta_1, \dots, \beta_l).
\end{eqnarray}
As a consequence, the following enumeration results hold.

\begin{prop}
Let $\mathbf{a}=(a_1, \dots, a_{k})$ , $\mathbf{b}=(b_1,
\dots, b_{l})$, where $a_j$, $b_j$  are  non-negative integers,
and from these arrays form partitions
$$
\mu=\langle 1^{a_1}\cdots k^{a_{k}}\rangle, \qquad
\nu=\langle 1^{b_1}\cdots
l^{b_l}\rangle
$$
Consider block matrices (\ref{abcd}) with constraints (\ref{CB})
and the further constraint that all diagonal entries of $A$ are
even (zero) while all those of $D$ are zero (even). Let
$N_{\mu \nu}^{\rm Sp}$ $(N_{\mu \nu}^{\rm O})$  denote the number of
such matrices with
$$
{\rm row}(A,B) = {\rm col}(A,C) = \mu, \qquad
{\rm row}(C,D) = {\rm col}(B,D) = \nu.
$$
We have
\begin{eqnarray*}
N_{\mu \nu}^{\rm Sp} & = & {\mathbb E}_{M \in
USp(2N)} \prod_{i=1}^{k} (\Sc_i(M))^{a_i}
 \prod_{j=1}^{l}
(\Rc_j(M))^{c_j} \\
N_{\mu \nu}^{\rm O} & = & {\mathbb E}_{M \in
O(N)} \prod_{i=1}^{k} (\Sc_i(M))^{a_i}
 \prod_{j=1}^{l}
(\Rc_j(M))^{c_j}.
\end{eqnarray*}
\end{prop}

\section{Relationship to non-intersecting paths} \label{path}
\setcounter{equation}{0}
The purpose of this section is to relate the random unitary matrix average
in (\ref{5.1}) and the unitary symplectic average in (\ref{lsp})
to configurations of weighted non-intersecting lattice paths. This then
allows us to give combinatorial interpretations to (\ref{1}) and (\ref{3})
relating to non-intersecting lattice paths rather than integer matrices.
That such interpretation are possible can be anticipated from the
RSK correspondence: the semi-standard tableaux therein have a well
known interpretation in terms of weighted non-intersecting lattice
paths (see e.g.~\cite{Sa}). However neither in the case of (\ref{5.1})
nor (\ref{lsp}) will our non-intersecting paths correspond to the
conventional ones related to semi-standard tableaux, and so a separate
discussion is warranted.

In relation to (\ref{5.1}), mark points on the $x$-axis at
$x=1,\dots,N$. Move each point to the line $y=1$ according to the rule
that each $x$ coordinate must either stay the same (weight unity) or
increase by one (weight $\alpha_1$), with the proviso that all $x$
coordinates must remain distinct. Connect the points between $y=0$ and
$y=1$ by segments, which must either be vertical (weight unity),
or right diagonal (weight $\alpha_1$). Repeat this procedure a total of
$n$ times, with each right diagonal segment at step $j$ weighted by
$\alpha_j$.

\begin{figure}[t]
\epsfxsize=7cm
\centerline{\epsfbox{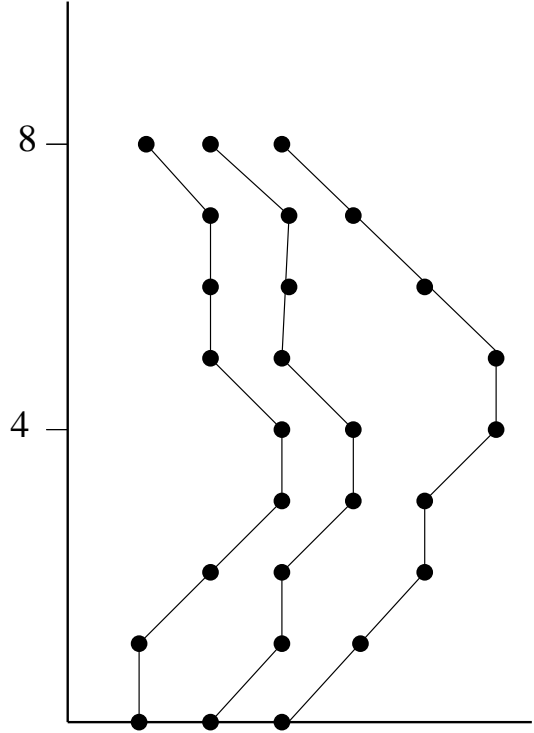}}
\caption{\label{F1} Example of non-intersecting lattice paths
corresponding to returning lock step vicious walkers. For the
first (last) $n$ steps the walkers must move to the right (left)
or stay stationary. The weight of the configuration shown here is
$\alpha_1^2 \alpha_2^2\alpha_3^2\alpha_4\beta_1^3\beta_2\beta_3\beta_4^2$.}
\end{figure}

After step $n$, perform another $n$ steps, but now with the segments
either vertical (weight unity) or left diagonal (weight $\beta_{2n+1-j}$
in step $n+j$). The segments again must not intersect, and are further
conditioned to return along $y=2n$ to the same $x$ coordinates
$(x=1,\dots,N)$ as they began (see Figure \ref{F1} for an example).
The resulting non-intersecting lattice paths are equivalent to a special
case of the lock-step model of vicious random walkers \cite{Fi84,Fo89}.
For general initial positions along $y=0$ $(l_1^{(0)},\dots, l_N^{(0)}$
say), and final positions along $y=2n$ $(l_1,\dots,l_N$ say), the
generating function $G_{2n}$ for the weighted paths can be written
as an $N \times N$ determinant according to
\begin{equation}
G_{2n}(l_1^{(0)},\dots,l_N^{(0)};l_1,\dots,l_N) =
\det \Big [ g_{2n}(l_j^{(0)};l_k) \Big ]_{j,k=1,\dots,N}
\end{equation}
where
\begin{equation}
g_{2n}(l^{(0)};l) = {1 \over 2 \pi} \int_{-\pi}^\pi
\prod_{j=1}^n(1 + \alpha_j e^{-i \theta_j})
(1 + \beta_j e^{i \theta_j}) e^{-i(l-l^{(0)})\theta} \, d \theta.
\end{equation}

{}From the well known identity
\begin{equation}
\det \Big [ {1 \over 2 \pi} \int_{-\pi}^\pi h(\theta) e^{-i(j-k) \theta}
\, d\theta \Big ]_{j,k=1,\dots,N} =
\Big \langle \prod_{j=1}^N h(\theta_j) \Big \rangle_{U(N)}
\end{equation}
we see that in the case of interest (initial positions $=$ final positions,
all one unit apart), we have
\begin{equation}
G_{2n}(l_1^{(0)},\dots,l_N^{(0)};l_1,\dots,l_N)
\Big |_{l_j^{(0)}=l_j=j \atop (j=1,\dots,N)}
= {\mathbb E}_{U \in U(N)} \prod_{j=1}^n \det (I_N + \alpha_j \bar{U})
\det (I_N + \beta_j U).
\end{equation}
The following interpretation of the average over $U(N)$ in (\ref{1}) is
now evident.

\begin{figure}[t]
\epsfxsize=7cm
\centerline{\epsfbox{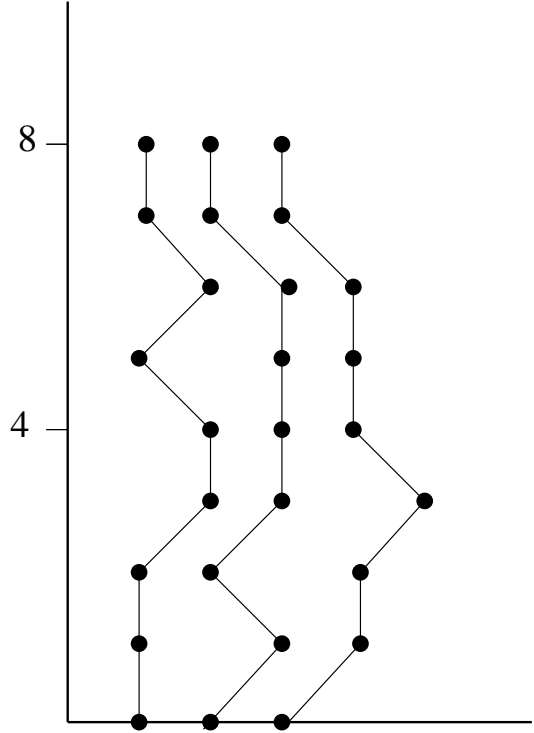}}
\caption{\label{F2} Example of non-intersecting lattice paths
with initial (final) spacings along $y=0$ $(y=2n)$ one unit apart
starting at $x=1$. All paths are restricted to $x\ge 1$, and segments
can be vertical or right diagonal for odd steps, and vertical or
left diagonal for even steps.
The weight of the configuration shown here is
$\alpha_1^3 \alpha_2^4\alpha_3^2\alpha_4^3$.}
\end{figure}

\begin{prop}\label{pt}
Consider the set of all non-intersecting paths of the type depicted in
Figure \ref{F1}. Impose the further constraint that the number of right
diagonal segments between $y=j-1$ and $y=j$ $(j=1,\dots,n)$ is equal to
$\mu_j$, and the total number of left diagonal segments between
$y=2n-j$ and $y=2n-j+1$ $(j=1,\dots,n)$ is equal to $\tilde{\mu}_j$.
The total number of such lattice paths is given by the RHS of
(\ref{1}).
\end{prop}

We note that in the setting of  Proposition \ref{pt} at most
$|\mu| = \sum_{j=1}^n \mu_j$ (which must equal $|\tilde{\mu}|$)
different walkers, counted to the left from $x=N$, can move.
Thus we see immediately that (\ref{1}) must be independent of $N$ for
$N \ge |\mu|$.

Let us now turn our attention to the non-intersecting lattice path
interpretation of (\ref{lsp}). For this mark points on the $x$-axis
at $x=1,\dots,N$. Also count steps $j=1,\dots,2n$ as the points
are moved to the lines $y=1,\dots,2n$ in order, and the corresponding
segments are drawn to connect the points. Let this process proceed by the
rule that the lattice paths must not intersect, and that at odd numbered
steps the segments must be either vertical or right diagonal, while at
even numbered steps the segments must either be vertical or left
diagonal. Weight the diagonal segments in steps $2j-1,2j$ by
$a_j, b_j$ respectively. It is further required that all points remain
to the right of the $y$-axis, which can be thought of as a wall
(see Figure \ref{F2} for an example). For general initial positions along
$y=0$ $(l_1^{(0)},\dots,l_N^{(0)})$ and finishing positions along
$y=2n$ $(l_1,\dots,l_N)$, all to the right of the wall, the generating
function $G_{2n}^{\rm wall}$ for the weighted paths has a
determinant form \cite{Fo89}. In the special case that
$\{a_i\} = \{b_i\}=\{\alpha_i\}$ (in any order), and that the initial
and final positions are spaced one unit apart starting at $x=1$,
this determinant can be expressed as the random matrix average
(\ref{lsp}) \cite{Fo01b}, thus providing the following interpretation of the
latter in terms of non-intersecting paths.

\begin{prop}
Consider the set of all non-intersecting lattice paths of the type depicted
in Figure \ref{F2}. Impose the further constraint that the sum of the
number of right diagonal segments at step $2j-1$ and the number of left
diagonal segments at step $2j$ is equal to $\mu_j$. The total number of
such paths is given by the RHS of (\ref{3}).
\end{prop}

\section{Relationship to plane partitions}\label{part}
\setcounter{equation}{0}
In this section
the random unitary matrix
average (\ref{e:charmom1}), the random unitary symplectic average
\eqref{lsp}, and the random orthogonal average
(\ref{2.16}) will be related
to the enumeration of certain classes of plane partitions
\cite{Mac, stanleyv2}.
A plane partition $\cal{P}$ is a finite
set of lattice points $\{(i, j,
k)\} \subseteq \natls^{3}$ with the property that if $(a, b, c)
\in \cal{P}$ and $1 \le i \le a$, $1 \le j \le b$, $1 \le k \le
c$, then $(i, j, k) \in \cal{P}$.  A plane partition is symmetric
if $(i, j, k) \in \cal{P}$ if and only if $(j, i, k) \in \cal{P}$.
The height of stack $(i, j)$ is the largest value of $k$ for which
there exists a point $(i, j, k)$ in the plane partition.

The study of plane partitions was initiated by MacMahon \cite{mm}
who proved that the generating function for plane partitions
fitting in the box \beq \label{mm1} \mathcal{B}(a, b, c) =\{(i, j,
k) | 1 \le i \le a, 1 \le j \le b, 1 \le k \le c \} \eeq is given by
\beq \label{mm2}
\prod_{i=1}^{a} \prod_{j=1}^{b} \prod_{k=1}^{c}
\frac{1-q^{i+j+k-1}}{1-q^{i+j+k-2}}.\eeq
Taking the limit $q \to 1$ gives for
the total
number of plane partitions fitting inside $\B(a, b, c)$, $\# {\cal P}(a,b,c)$
say, the evaluation
\begin{equation}\label{6.3}
\# {\cal P}(a,b,c) = \prod_{i=1}^{a}
\prod_{j=1}^{b} \prod_{k=1}^{c} \frac{i+j+k-1}{i+j+k-2}.
\end{equation}

We can express the generating function (\ref{mm2}) for plane partitions
in terms of Schur functions.  Thus from the
combinatorial definition of  Schur functions $s_\lambda(x_1,\dots,x_n)$
as a sum over weighted semi-standard
tableaux of shape $\lambda$ and content
$\{1,\dots,n\}$, it follows that
$$
s_{b^a}(q^{a+c}, q^{a+c-1}, \dots, q)
$$
where
$b^a$ denotes a partition with $a$ parts all of which are equal to
$b$, is the generating function for plane partitions strictly
decreasing down columns with exactly $a$ rows each of length $b$
and with the largest stack of height less than or equal to $a+c$.
The strictly decreasing constraint can be eased by
removing $a-i+1$ from the boxes in row $i$ to bijectively obtain a plane
partition which is a subset of $\mathcal{B}(a, b, c)$.
Consequently an alternative expression for \eqref{mm2} is given by
\beq \label{mm3} q^{-ba(a+1)/2} s_{b^a}(q^{a+c}, q^{a+c-1}, \dots,
q).\eeq
Making use of (\ref{e:charmom1}) shows that for
$a=n$, $c=m$ and
$b=N$, \eqref{mm3} can be expressed as
$$
{\mathbb E}_{U \in U(N)}
\prod_{j=1}^m \det(I_N + \alpha_j U) \prod_{j=1}^n
\det(I_N + \overline{
U}/\beta_j
)
$$
with $\alpha_i=q^{n+i}$ and $\beta_i=q^{i}$.  Taking the limit $q \to 1$ gives
the sought relationship between (\ref{6.3}) and a random matrix average.

\begin{prop}
Denote by $P_M^*(z)$ the (reciprocal) characteristic polynomial of
$M \in  U(N)$, so that $P_M(z) = \det (I_N - z M)$.
For $|z|=1$ we have
\begin{equation}\label{mz}
\# {\cal P}(a,b,c)
=  {\mathbb E}_{U \in U(c)} \Big ( P_M^*(z) \Big )^b
\Big ( \overline{P_M^*(z)} \Big )^a.
\end{equation}
\end{prop}

\noindent
We remark that in \cite{St01}, the case $a=b$ of the RHS of
(\ref{mz}) has been given a combinatorial interpretation involving
two-rowed lexicographic arrays.

The random matrix average (\ref{mz}) is an example of a class
of multi-dimensional integrals possessing gamma function
evaluations (see e.g.~\cite{Fo02}). This implies
\begin{eqnarray}\label{bG}
\# {\cal P}(a,b,c) & = &
\prod_{j=0}^{c-1} {\Gamma(a+b+1+j) \Gamma(2+j) \over
\Gamma(a+1+j) \Gamma(b+1+j) } \nonumber \\
& = &
{G(1+a+b+c) \over G(1+a+b)}
{G(1+a) \over G(1+a+c)} {G(1+b) \over G(1+b+c)} G(c+2)
\end{eqnarray}
where $G(z)$ is the Barnes $G$-function, related to the gamma function
by the functional equation $G(z+1) = \Gamma(z) G(z)$, $G(1)=1$.
In contrast to the formula (\ref{6.3}), the latter formula is
well suited for asymptotic analysis.

Consider now symmetric plane partitions fitting inside the box
(\ref{mm1}) with $b=a$. The generating function for such plane
partitions, with the additional constraint that the heights of all
stacks on the diagonal are even and bounded by $2c$, has the
product form \cite{stanleyv2}
\begin{equation}\label{qs}
\prod_{1 \le i \le j \le a} {1 - q^{i+j+2c} \over 1 - q^{i + j} }.
\end{equation}
Here only points on and above the diagonal are weighted.
Denoting their total number by $\# {\cal P}^{\rm sym}_{\rm e}(a,2c)$, we see by
taking the limit $q \to 1$ that \cite{dsv}
\begin{equation}\label{qs1}
\# {\cal P}^{\rm sym}_{\rm e}(a,2c) = \prod_{1 \le i \le j \le a}
{i+j+2c \over i + j}.
\end{equation}

The generating function (\ref{qs}) and thus counting formula (\ref{qs1}) can
be expressed in terms of Schur functions by making use of a bijection
between tableaux and symmetrical plane partitions. Consider then
semi-standard tableaux of content $\{1,\dots,a\}$. Suppose furthermore
that each
row is of even length, and the first row is constrained to be less than
or equal to $2c$. From this construct the diagonal and upper triangular
portion of a symmetrical plane partition by associating with
grid points $(i,j)$, $i \le j$,
stacks of height
$$
h_{i,j} = \#i{\rm '}s + \#(i+1){\rm '}s + \cdots +
\#(a+i-j){\rm '}s \: \: {\rm in} \: \:
{\rm row} \: \: i \: \: {\rm of} \: \: {\rm the} \: \:
{\rm tableaux}.
$$
By weighting each square labelled $j$ in the
tableaux by $q^j$, we see that each stack at grid point $(i,j)$ is weighted
$q^{h_{i,j}}$, and furthermore on the diagonal this weight is $q^{\lambda_i}$
where $\lambda_i$ is the length of row $i$ of the tableau. It follows
that (\ref{qs}) can be expressed in terms of Schur functions according to
\begin{equation}\label{5.8a}
\sum_{ \substack{\lambda \subseteq (2c )^a \\ \lambda \,
\text{even}}} s_{\lambda}(q^a, q^{a-1},\dots,q) .
\end{equation}
Recalling (\ref{lsp}) gives the sought relationship with an
average over the unitary symplectic group.

\begin{prop}
Denote by $P_M^{\rm Sp}(z)$ the characteristic polynomial of
$M \in USp(2c)$ so that $P_M^{\rm Sp}(z) = \det (\lambda I_{2c} -
M)$. We have
$$
\# {\cal P}^{\rm sym}_{\rm e}(a,2c) = {\mathbb E}_{M \in
USp(2c)}\Big (P_M^{\rm Sp}(-1) \Big )^a.
$$
\end{prop}

Analogous to (\ref{bG}), the above average is of a type which
admits a gamma function evaluation (see e.g.~\cite{Fo02}). This allows
(\ref{qs1}) to be written
\begin{eqnarray*}
\# {\cal P}^{\rm sym}_{\rm e}(a,2c) & = &
2^{2ca}\prod_{j=1}^{c}\frac{\Gamma(1+c+j)
\Gamma(\frac{1}{2}+a+j)}{\Gamma(1+c+a+j) \Gamma(\frac{1}{2}+j)} \\
 & = &
2^{2ca} {G(2+2c) \over G(2+c)}
{G(2+c+a) \over G(2+2c+a)} {G({3 \over 2} + a + c) \over
G({3 \over 2} + a) } {G({3 \over 2}) \over G({3 \over 2}+c)}.
\end{eqnarray*}

Finally, consider the same weighted symmetric plane partitions as
described in the paragraph below (\ref{qs1}), but with stacks now
restricted to heights $\le c$ and without the constraint that the
stacks on the diagonal be even. The generating function is
$$
\sum_{\lambda \subseteq c^a} s_\lambda(q^a, q^{a-1},\dots,q).
$$
In an obvious notation, after taking the limit $q \to 1$ and recalling
(\ref{7.1b}) we obtain
$$
\# {\cal P}^{\rm sym}(a,c) = {\mathbb E}_{M \in O(c)}\Big
(P_M^{\rm O * }(-1) \Big )^{a+1}.
$$

\section{Counting formulas associated with the
characteristic polynomials  of random
Hermitian and Wishart matrices} \label{herm}

\setcounter{equation}{0} Let $X$ be an $n \times n$ Hermitian
matrix. Let the diagonal elements $x_{ii}$ be chosen independently
according to a probability distribution ${\mathcal D}_1$, with the
property that
\begin{equation}\label{a5.1}
{\mathbb E}_{{\mathcal D}_1} x_{ii} = 0.
\end{equation}
Let the upper triangular elements $x_{ij}$, $i < j$, (which may be
complex) be chosen independently according to a probability
distribution ${\mathcal D}_2$, with the properties that
\begin{equation}\label{a5.2}
{\mathbb E}_{{\mathcal D}_2} x_{ij} = 0, \qquad {\mathbb
E}_{{\mathcal D}_2} |x_{ij}|^2 = \sigma_2^2.
\end{equation}
The following result, generalizing a number of results presented
in \cite{DiGa}\footnote{We take this opportunity to correct the
statement of Theorem 15 in \cite{DiGa}.  It should read:
$\E_{\mu_N}(P_M^{2k}(x))= h_N^{(k)}(x), $ where $h_N^{(k)}$ are
orthogonal polynomials with respect to the weight
$(t-x)^{2k}e^{-t^2}$.}, gives the expected value of the
characteristic polynomial of $X$ in terms of the classical Hermite
polynomial $H_n(x)$.

\begin{prop}\label{p6}
Let $X$ be Hermitian and specified in terms of ${\mathcal D}_1$
and ${\mathcal D}_2$ as above. With the monic rescaled Hermite
polynomial specified in terms of the classical Hermite polynomial
by
$$
h_n(x) = 2^{-n/2} H_n \Big ( {x \over \sqrt{2} } \Big )
$$
we have
\begin{equation}
{\mathbb E}_{{\mathcal D}_1,{\mathcal D}_2} \det ( \lambda I_N -
X) = \sigma_2^N h_N \Big ( {\lambda \over \sigma_2} \Big ).
\end{equation}
\end{prop}

\noindent
Proof. \quad By definition
\begin{equation}\label{a6.1}
\det ( \lambda I_N - X) = \sum_{P \in S_N} \varepsilon(P)
\prod_{l=1}^N ( \lambda_{l,P(l)} - x_{l,P(l)} ),
\end{equation}
where $\varepsilon(P)$ denotes the parity of $P$ and
$$\lambda_{i, j}=\begin{cases}
\lambda  &\text{if $i=j$;}\\
0, &\text{if $i \ne j$}.
\end{cases}$$

The specifications (\ref{a5.1}), (\ref{a5.2}) tell us that the
only non-zero terms in (\ref{a6.1}) after averaging over
${\mathcal D}_1,{\mathcal D}_2$ are those for which $P$ consists
entirely of fixed points $(P(j)=j)$ and 2-cycles $(P(j_1) = j_2$
and $P(j_2) = j_1$, $j_1 \ne j_2$). Let there then be $N - 2j$
fixed points and $j$ 2-cycles. Such permutations have parity
$(-1)^j$. Each fixed point contributes a factor $\lambda$, while
each 2-cycle contributes $\sigma_2^2$. As the number of ways of
choosing the $(N-2j)$ fixed points and the $j$ 2-cycles is
$$
\Big ( {N \atop 2j} \Big ) {(2j)! \over 2^j j!}
$$
we see that
$$
{\mathbb E}_{{\mathcal D}_1,{\mathcal D}_2} \det ( \lambda I_N -
X) = \sum_{j=0}^{[N/2]} (-1)^j \Big ( {N \atop 2j} \Big ) {(2j)!
\over 2^j j!}
 \lambda^{N-2j} \sigma_2^{2j}.
$$
But this is precisely the power series expansion of the polynomial
in question.
\hfill $\square$

\medskip
Next we turn our attention to the mean characteristic polynomial
of so called chiral matrices, that is matrices of the form
\begin{equation}\label{K}
K := \left [ \begin{array}{cc} 0_{n \times n} & X_{n \times p}
\\ (X^\dagger)_{p \times n} & 0_{p \times p} \end{array}
\right ],
\end{equation}
where we require that $n \ge p$. These matrices have exactly $n-p$
zero eigenvalues, with the remaining $2p$ eigenvalues given by
$\pm$ the positive square roots of the eigenvalues of the
non-negative matrix $X^\dagger X$. We specify that the elements of
$X$ be chosen independently according to a probability
distribution ${\mathcal D}$ with the properties
\begin{equation}\label{f1}
{\mathbb E}_{\mathcal D} x_{ij} = 0, \qquad {\mathbb E}_{\mathcal
D} |x_{ij} |^2 = \sigma^2.
\end{equation}
The analogue of Proposition \ref{p6} can readily be deduced.

\begin{prop}
Let $K$ be a random chiral matrix as specified above. We have
\begin{equation}\label{f2}
{\mathbb E}_{\mathcal D} \det (\lambda I_{n+p} - K) = p!
\sigma^{2p} \lambda^{n-p} L_p^{n-p}((\lambda/\sigma)^2),
\end{equation}
where $L_m^a(x)$ denotes the classical Laguerre polynomial.
\end{prop}

\noindent
Proof. \quad Because of the first specification in (\ref{f1}), we see that
in the analogue of (\ref{a6.1}) for $\det (\lambda I_{n+p} - K)$, the
only terms after averaging will again result entirely from fixed points
and 2-cycles. The fact that $K$ has $n-p$ zero eigenvalues implies we
require there be a minimum of $n-p$ fixed points. Thus we must consider
the cases that the number of fixed points is equal to $n+p-2j$ and
the number of 2-cycles is equal to $j$ for each $j=0,\dots,p$.

Not all permutations with these specifications give a non-zero
contribution. For the latter, because of the zero blocks in
(\ref{K}), in relation to the 2-cycles $(j_1 j_2)$
we require
$$
j_1 \in \{1,\dots,n\}, \qquad j_2 \in \{n+1,\dots,n+p\}.
$$
Hence we must choose $(n-p)+(p-j)$ fixed points from the first of these
sets and $(p-j)$ from the second. The number of distinct ways to do this
is
$$
\Big ( {n \atop n-j} \Big ) \Big ( {p \atop p-j} \Big ).
$$
The $j$ 2-cycles can then be chosen in $j!$ different ways. Since again
each fixed point contributes a factor $\lambda$, while each 2-cycle contributes
a factor $\sigma^2$, we see that
\begin{eqnarray*}
{\mathbb E}_{\mathcal D} \det (\lambda I_{n+p} - K) & = &
\sigma^{2p} \lambda^{n-p} \sum_{j=0}^p \Big ( {n \atop n-j} \Big )
\Big ( {p \atop p-j} \Big ) j!
\Big ( {\lambda \over \sigma} \Big )^{2(p-j)} \nonumber \\
& = &
p! \sigma^{2p} \lambda^{n-p} \sum_{j=0}^p
\Big ( {n \atop p-j} \Big ) {1 \over j!}
\Big ( {\lambda \over \sigma} \Big )^{2j}.
\end{eqnarray*}
The sum in this expression is precisely $L_p^{n-p}((\lambda/\sigma)^2)$.
\hfill $\square$

\medskip
Due to the relationship between the chiral matrices $X^\dagger X$ as noted
below (\ref{K}), we have the following result for the expected value
of the characteristic polynomial for Wishart matrices.
\begin{cor}
Let $X$ be a $n \times p$ matrix with elements independently
chosen according to the distribution ${\mathcal D}$ with the
properties (\ref{f1}). We have
\begin{equation}
{\mathbb E}_{\mathcal D} \det (\lambda I_p - X^\dagger X) = p!
\sigma^p L_p^{n-p}(\lambda/\sigma).
\end{equation}
\end{cor}

\section*{Acknowledgements}
We thank the organizers of the program `Random matrix methods in
number theory', held at the Newton Institute during the first half
of 2004, for inviting us to participate and so facilitating our
collaboration. The work of PJF was supported by the Australian
Research Council.  The work of AG was supported in part by NSF
Postdoctoral Fellowship.

\end{document}